\documentclass[11pt]{article}

\usepackage[margin = 1.0in]{geometry}   
\usepackage{geometry}  			
\usepackage{amssymb}
\usepackage{amsmath}
\usepackage{palatino}
\usepackage{graphicx}
\usepackage[english]{babel}
\usepackage{amscd}
\usepackage{amsthm}
\usepackage{color}
\usepackage{enumerate}
\usepackage{hyperref}
\usepackage{subcaption}
\usepackage{comment}

\newtheorem{theorem}{Theorem}

\newtheorem{corollary}[theorem]{Corollary}

\newtheorem{proposition}{Proposition}

\newtheorem{remark}{Remark}

\newcommand{\Prb}{\mathbb{P}}
\newcommand{\Exp}{\mathbb{E}}

\newcommand{\Var}{\operatorname{Var}}
\newcommand{\Hyp}{\operatorname{Hyp}}

\newcommand{\Bin}{\operatorname{Bin}}
\newcommand{\Pois}{\operatorname{Pois}}

\newcommand{\bp}{\mathbf{p}}
\newcommand{\bbZ}{\mathbb{Z}}
\newcommand{\veps}{\varepsilon}

\title{On the probability that a discrete IFR random variable exceeds its mean}

\author{
Jianhang Ai\thanks{Faculty of Electrical Engineering, Czech Technical University, Karlovo N\'am\v{e}st\'i 13, 12135, Prague 2, Czech Republic, e-mail: ai.jianhang@fel.cvut.cz}
\and
Vassilis Karagiannis\thanks{Aristotle University of Thessaloniki, 
Department of Mathematics, 
541 24 Thessaloniki, Greece, e-mail: vkdstat@math.auth.gr}
\and
Ond\v{r}ej Ku\v{z}elka\thanks{Faculty of Electrical Engineering, Czech Technical University, Karlovo N\'am\v{e}st\'i 13, 12135, Prague 2, Czech Republic, e-mail: ondrej.kuzelka@fel.cvut.cz}
\and
Christos Pelekis\thanks{Aristotle University of Thessaloniki, 
Department of Mathematics, 
541 24 Thessaloniki, Greece, e-mail:  cpelekis@math.auth.gr }
}

\begin{document}

\maketitle

\begin{abstract}  
Let $X$ be an integer-valued random variable with $\mathbb{E}(|X|)<\infty$ whose survival probabilities $S_j=\mathbb{P}(X\ge j), j\in\mathbb{Z}$, form a log-concave sequence. 
Suppose  that for some  integer $m$ and some $\gamma\in (0, 1/2]$ it holds 
$\mathbb{E}(X)\ge m$ and $\mathbb{P}(X=m+1)\ge\gamma\cdot\mathbb{P}(X=m)> 0$. It is shown that  
\[
\mathbb{P}(X \ge m+1) \,\ge\, \frac{\gamma}{1+2\gamma} \, ,
\]
with the bound  being best possible. 
As a consequence, it is shown that if $H \sim \text{Hyp}(n, i, k)$ is a hypergeometric random variable for which both $ik/n$ and $(n-i)(n-k)/n$ 
are at least $\log(4/3)$, then it holds $\mathbb{P}(H >\mathbb{E}(H)) \ge 1/4$ when $\mathbb{E}(H)$ is not an integer, and  $\mathbb{P}(H >\mathbb{E}(H)) \ge 1/6$ when $\mathbb{E}(H)$ is an integer, both bounds being best possible. 
\end{abstract}

\vspace{0.2cm}


\section{Introduction}
\label{intro}

We shall be interested in constant lower bounds on the probability that an 
integer-valued random variable exceeds its mean.
Without additional ``regularity assumptions" on the distribution, no positive
universal lower bound is possible. For instance, let $n\ge2$ and consider the random variable $X$ that takes values in the set $\{0,n\}$ with probabilities $1-1/n$ and $1/n$, respectively. 
Then $\Exp(X)=1$, whereas
$\Prb(X>\Exp(X))=1/n\to0$ as $n\to\infty$.

This motivates the search for structural conditions on the probability distribution under which
the probability of exceeding the mean is bounded away from zero.
In this article, we establish a sharp quantitative lower bound
on the aforementioned  probability under an increasing failure rate (IFR) assumption
and a ``local comparison" between two adjacent probability masses
near the mean. 
Before presenting our main result, let us recall some basic facts and definitions from the class of IFR random variables that will be used in the sequel.

\subsection{Notation and definitions}
 
Throughout the text, given a real number $x$, we write $\lfloor x\rfloor = \sup\{m\in\mathbb{Z} : m\le x\}$ and  $x_+ = \max\{0,x\}$. 
If $X$ and $Y$ are two random variables, we write $X\sim Y$ when they have the same distribution. 
For integers $a\le b$, we write  $[a..b]=\{a,a+1,\ldots,b\}$, abbreviate $[1..n]$ to $[n]$, and refer to $[a..b]$ as a \emph{discrete interval}. 
If $a=0$ and $b=\infty$, the corresponding discrete interval is  denoted $\mathbb{Z}_{\ge 0}$, and 
if $a=-\infty$ and $b=\infty$, the corresponding interval is denoted $\bbZ$.

The \emph{support} of a non-negative sequence of real numbers $\mathbf{p}=\{p_j\}_{j\in\bbZ}$ is the set $\text{supp}(\mathbf{p}) = \{j\in \bbZ : p_j >0\}$. The  sequence $\mathbf{p}=\{p_j\}_{j\in\bbZ}$ is called \emph{log-concave} if its  support is a discrete interval and it satisfies 
\[
 p_j^2\,\ge\, p_{j-1}p_{j+1}\, ,\,\,\text{ for all } \, \, j\in\text{supp}(\mathbf{p}) \, .
\]
Equivalently, the sequence $\{p_{j+1}/p_j\}_{j\in\text{supp}(p)}$ of the successive ratios is non-increasing.
A random variable $X$ taking values in the set of non-negative integers is referred to as a \emph{$\bbZ_{\ge 0}$-valued} random variable.  
An integer-valued random variable is also referred to as $\bbZ$-valued. A $\bbZ$-valued  random variable is 
called \emph{log-concave} if its probability distribution, denoted $\bp_X=\{p_j\}_{j\in\bbZ}$, is a log-concave sequence. 
The \emph{support} of $X$ is the support of $\bp_X$. 
The \emph{survival sequence} of $X$ is the sequence $\mathcal{S}_X=\{S_j\}_{j\in\bbZ}$, where $S_j = \Prb(X\ge j)$, for $j\in\bbZ$. 
A $\bbZ$-valued random variable will be referred to as  \emph{IFR (increasing failure rate)} if its 
survival sequence form a log-concave sequence  
(see~\cite{Shaked_Shantikumar}). It is a well-known fact (see, for example,~\cite[Lemma~4.1]{Bentkus}), and rather straightforward to verify, that if $X$ is log-concave then it is IFR. 

\subsection{Related work and main results}

Let $X$ be an integer-valued random variable with finite mean. A natural question is to ask how small the probability that $X$ exceeds its mean can be. 
What is a sharp lower bound on the tail probability $\Prb(X > \Exp(X))$? 

As already mentioned, some additional "regularity assumptions" on the probability distribution are needed in order to hope for non-trivial lower bounds on the above-mentioned tail probability. 

In the continuous setting, a well-studied class of  "regular" probability distributions  are the log-concave distributions.
Recall that a continuous real-valued random variable $X$,  with density $f$, is called log-concave if 
\[
f(\alpha x + (1 - \alpha)y) \ge  f(x)^{\alpha}\cdot f(y)^{1-\alpha} \, ,
\]
holds true for $x,y$ in the support of $X$ and all $\alpha\in [0,1]$. 
A classical result going back to  Gr\"unbaum \cite{Grunbaum} (see also~\cite[Lemma~5.4]{Lovasz_Vempala}) asserts that for a continuous and log-concave real-valued random variable $X$ it holds 
\[
    \Prb(X > \Exp(X)) \,\ge\, \frac{1}{e} \, ,
\]
with the constant $1/e$ being sharp and with the bound being attained by an exponential random variable of parameter $1$. By reflection, this bound is complemented by the bound 
$\Prb(X< \Exp(X) )\geq \frac1e$, and both bounds 
illustrate how log-concavity prevents the tails  from having  small probability. In the case of non-negative, continuous random variables,  a refinement of Gr\"unbaum's bound holds true for the more general class of IFR continuous distributions which, by definition, are continuous random variables $X$ having log-concave survival function  $S(t)=\Prb(X> t)$.  
In particular, 
a classical  result of  Barlow, Proschan and Marshall (see~\cite[Theorem~3.8]{Barlow_Marshall} and~\cite[Theorem~4.4]{Barlow_Proschan}) asserts that 
for a nonnegative-real-valued and IFR continuous random variable $X$ it holds 
\begin{equation}\label{BPM}
\Prb(X > t) \ge e^{-\frac{t}{\Exp(X)}} \, , \,\, \text{ for } \,\, 0\le t < \Exp(X) \, .
\end{equation}
Hence, letting $t\to\Exp(X)$ yields the bound 
$\Prb(X>\Exp(X))\ge \frac{1}{e}$ and, as in the log-concave  setting, the IFR  assumption on a continuous  distribution guarantees a universal lower bound on the probability of exceeding the mean.

In the discrete setting, it has been shown by   
Shah (see~\cite[Theorem~2]{Shah}) that if $X$ is  a $\bbZ_{\ge 0}$-valued and IFR random variable with  mean $\mu= \Exp(X)>0$ then for every integer $0\le k \le \mu$ it holds 
\begin{equation}\label{eq:Shah_bound}
  \Prb\big(X\,\ge\, k \big) \,\ge\,  \left( \frac{\mu}{1+\mu}\right)^k \,\ge\, e^{-k/\mu} \, .   
\end{equation}
In particular, a $\bbZ_{\ge 0}$-valued and IFR random variable satisfies  $\Prb\big(X > \Exp(X)-1\big) \ge \frac{1}{e}$.
The first bound~\eqref{eq:Shah_bound}  is sharp, and is attained by a geometric random variable $X$ with mass 
$\Prb(X=j)= \frac{1}{1+\mu} \cdot \left( \frac{\mu}{1+\mu}\right)^j$, for $j\ge 0$. 
In the case of integer-valued random variables, the following  bound  has been reported  in~\cite{Alqasem_et_al}.

\begin{theorem}[Alqasem et al.~\cite{Alqasem_et_al}]
\label{thm:Alqasem_et_al}
Let $X$ be a $\bbZ$-valued and log-concave random variable. Then 
\[
 \Prb\big(X \,>\, \Exp(X)-1\big) \,\ge\, \frac{1}{e} \, .
\]
The constant $1/e$ is best possible. 
\end{theorem}

Note that both~\eqref{eq:Shah_bound} and the bound provided by Theorem~\ref{thm:Alqasem_et_al} concern lower bounds on the probability that the random variable under consideration takes a value that is above a fixed value below the mean.

In this article we complement  the aforementioned results by establishing a constant lower bound on the probability that an  integer-valued random variable takes a value "right above" its mean, under an additional "regularity condition". 
Such a condition is necessary even within the class of log-concave distributions. Indeed, let $0<\veps <1/3$ and consider the random variable $X$ taking values in the set $\{0,1,2\}$ with probabilities $\veps, 1-2\veps, \veps$, respectively. 
Then $X$ is, on the one hand, log-concave and so IFR but, on the other hand, it  
satisfies $\Exp(X)=1$ and   $\Prb(X > 1)=\veps\to 0$.

Our main result quantifies how a lower bound on the ratio
of two consecutive probability masses near the mean
yields a lower bound on the corresponding upper tail, and reads as follows.

\begin{theorem}\label{thm:main1}
Let $X$ be a $\bbZ$-valued IFR random variable such that 
$\Exp(|X|)<\infty$, with  probability distribution  $\bp_X=\{p_j\}_{j}$.
Assume further that there exists a real number  $\gamma\in(0,1/2]$ and an integer $m$ such that 
\begin{equation}\label{eq:criterion-assumptions}
\Exp(X)\ge m
\qquad\text{and}\qquad
p_{m+1}\ge\gamma\cdot p_m>0 \, .
\end{equation}
Then
\begin{equation}\label{eq:quantitative_ifr}
\Prb(X\ge m+1)\ge\frac{\gamma}{1+2\gamma}.
\end{equation}
The bound is sharp for every $\gamma\in(0,1/2]$. 
In particular, if
\begin{equation}\label{eq:left_skewed}
\Prb\bigl(X=\lfloor\Exp(X)\rfloor+1\bigr)
\ge
\frac12\cdot \Prb\bigl(X=\lfloor\Exp(X)\rfloor\bigr) \, ,
\end{equation}
then
\[
\Prb(X>\Exp(X))\ge\frac14 \, .
\]
\end{theorem}

For each $\gamma\in(0,1/2]$, equality
in~\eqref{eq:quantitative_ifr} is attained by the
log-concave distribution on $\{m-1,m,m+1\}$ with masses
$\frac{\gamma}{1+2\gamma},\,\frac{1}{1+2\gamma},\,\frac{\gamma}{1+2\gamma}$, 
respectively, whose  expected value is equal to $m$, and which satisfies 
$p_{m+1}/p_m=\gamma$.

\subsection{The hypergeometric distribution}

Theorem~\ref{thm:main1} applies to a variety of discrete probability distributions, with a notable case being the hypergeometric. 
Fix a positive integer $n\ge 2$, and $i,k\in [n]$, and  consider an urn consisting of $i$ black and $n-i$ white marbles. 
Let $H$ denote the number of black marbles in a sample without replacement of size $k$ from the urn. 
Then $H$ is  a hypergeometric random variable, denoted $H\sim\Hyp(n,i,k)$, having probability mass 
$\Prb(H=j)=\frac{\binom{i}{j}\binom{n-i}{k-j}}{\binom{n}{k}}$, for  $j\in [\max\{0,k-(n-i)\}\,..\,\min\{k,i\}]$. 
Recall that $\Exp(H)= \frac{ik}{n}$ and $\text{Var}(H) = k\cdot\frac{i}{n}\cdot\frac{n-i}{n}\cdot\frac{n-k}{n-1}$.

By a result of Vatutin and Mikhailov (\cite{Vatutin_Mikhailov}, see also~\cite{Hui_Park}), a hypergeometric random variable is a sum of independent Bernoulli random variables. Hence the hypergeometric distribution is amenable to Shevtsova's (see~\cite{Shevtsova}) Berry-Essesn-type inequality  for non-identically distributed summands.
More concretely, let $F$ denote the distribution function of $H$, where $\sigma=\sqrt{\Var(H)}>0$, and set $G(x)=\Phi\bigl((x-\Exp(H))/\sigma\bigr)$, where $\Phi$ is the distribution function of the standard normal random variable.
Mattner and Schulz~\cite[Formula~(17)]{Mattner_Schulz} observed that Shevtsova's result implies  $\sup_{x\in\mathbb{R}}|F(x)-G(x)|<0.5583/\sigma$. More recently, 
Mattner~\cite[Theorem~1.1 and Remark~1.2(b), Formula~(4)]{Mattner} improved the constant $0.5583$ to $0.3440$ in the aforementioned bound, by applying a structural refinement of Shevtsova's inequality to Bernoulli summands. In particular, the following holds:
\begin{equation}\label{eq:MS}
\Prb\big(H\,>\, \Exp(H)\big) \,\ge\, \frac{1}{2} - \frac{0.3440}{\sqrt{\Var(H)}} \, .
\end{equation}
Consequently, $\Prb(H>\Exp(H))\ge\frac14$ whenever $\Var(H)\ge(4\cdot0.3440)^2=1.893376$. 
However, as has been suggested in~\cite{Ai_Pelekis},  there is space  for improvements upon the bound provided by~\eqref{eq:MS}, which does not cover the case of "small variance". 
Using our main result, we complement this bound by establishing the following.

\begin{theorem}\label{thm:main_hyp}
Let $n\ge 2$ be an integer and fix $i,k\in \{1,\ldots,n\}$. Let  $H\sim\Hyp(n,i,k)$ be such that 
$\frac{ik}{n} \ge \log(4/3)$ and $\frac{(n-i)(n-k)}{n} \ge \log(4/3)$. 
\begin{enumerate}
\item If\, $\Exp(H)$ is \emph{not} an integer, then 
$\Prb\big(H\,>\,\Exp(H) \big)\,\ge\, \frac14$. 
\item If\, $\Exp(H)$ is an integer, then 
$\Prb\big(H\,>\,\Exp(H) \big)\,\ge\, \frac{1}{4 + 2/\Exp(H)} \,\ge\, \frac16$. 
\end{enumerate}
All constants $1/4, 1/6$ and $\log(4/3)$ are best possible. 
\end{theorem}

\subsection{Organization}

The remaining part of our article is organized as follows.
In Section~\ref{sec:main1} we prove our main result
(Theorem~\ref{thm:main1}) using basic properties of the log-concavity of the
survival sequence.
Section~\ref{sec:hypergeometric} establishes the hypergeometric
bounds of Theorem~\ref{thm:main_hyp} and the optimality of
their constants.
We conclude with Section~\ref{sec:applications}, in which
we apply our main result to several discrete probability distributions, including Poisson-binomial distributions,  differences of independent binomial random variables, 
as well as Poisson, negative binomial
and negative hypergeometric distributions.

\section{Proof of Theorem~\ref{thm:main1}} 
\label{sec:main1}

For every $j\in\bbZ$, set $S_j=\Prb(X\ge j)$.
Since $X$ is IFR,  the sequence $\{S_j\}_{j\in\bbZ}$ is log-concave. We aim to show that
$S_{m+1}\ge\frac{\gamma}{1+2\gamma}$. 
Set $x=S_{m+1}/p_m$ and note that,  since $S_{m+1}\ge p_{m+1}$,
condition~\eqref{eq:criterion-assumptions} implies  $x =S_{m+1}/p_m\ge p_{m+1}/p_m \ge\gamma$.
We distinguish two cases.

Assume first that $x\le1$; hence $x\in[\gamma,1]$.
We claim that it holds 
\begin{equation*}
  \Exp\left((X-m)_+ \right)=  \sum_{j\ge m+1}S_j \,\le\, \frac{S_{m+1}^2}{p_{m+1}} \, .
\end{equation*}
To prove the claim, note that if $p_{m+2}=0$, then it holds $S_{m+1}=p_{m+1}$ and $\sum_{j\ge m+1}S_j = S_{m+1}= \frac{S_{m+1}^2}{p_{m+1}}$. So we may assume that $m+2$ belongs to the support of $p_X$. In this case, if $j\ge m+1$ is such that $S_j>0$, 
the log-concavity of $\{S_j\}_j$ implies that 
\[
S_j = S_{m+1}\cdot\prod_{i=1}^{j-m-1} \frac{S_{m+1+i}}{S_{m+i}} \le S_{m+1}\cdot \left(\frac{S_{m+2}}{S_{m+1}} \right)^{j-m-1} \, ,
\]
and therefore we have    
\begin{equation*}
  \sum_{j\ge m+1} S_j\,\le\, S_{m+1}\cdot  \sum_{j\ge m+1}\left(1-\frac{p_{m+1}}{S_{m+1}} \right)^{j-m-1} \,=\, 
 \frac{S_{m+1}^2}{p_{m+1}} \, . 
\end{equation*}
The claim follows. 
Now combine the claim with the estimate  
\[
 \Exp\left((m-X)_+\right)\ge\Prb(X\le m-1)=1-p_m-S_{m+1} \, ,
\]
to conclude that 
\begin{eqnarray*}
0&\le& \Exp(X-m) \\
&=& \Exp\left((X-m)_+ \right) - \Exp\left((m-X)_+ \right)\\
&\le&  \frac{S_{m+1}^2}{p_{m+1}}+p_m+S_{m+1}-1\\
&\le& p_m\left(\frac{x^2}{\gamma}+1+x\right)-1\\
&\le& \frac{1+2\gamma}{\gamma}\cdot S_{m+1}-1 \, ,
\end{eqnarray*}
where the last estimate follows from the inequality 
$\frac{x^2}{\gamma}+1+x-\frac{1+2\gamma}{\gamma}x=\frac{(x-\gamma)(x-1)}{\gamma}
\le0$, valid for $x\in[\gamma,1]$. 
We have shown that $S_{m+1}\ge\gamma/(1+2\gamma)$, and the result follows.

So we may assume from now on that  $x=\frac{S_{m+1}}{p_m}>1$.  
To simplify notation, set  $d=\lfloor \frac{S_{m+1}}{p_m}\rfloor\ge 1$.   
As before, the log-concavity of $\{S_j\}_j$ implies that for every $j\ge m+1$ for which $S_j>0$ it holds 
$S_j\le S_{m+1}\cdot\left(\frac{S_{m+1}}{S_{m}} \right)^{j-m-1}$. Similarly, if $j\le m$ we have 
\[
S_j = \frac{S_{m+1}}{\prod_{i=j}^{m} \frac{S_{i+1}}{S_{i}} }  \le \frac{S_{m+1}}{\left(\frac{S_{m+1}}{S_{m}}\right)^{m-j+1}} =
S_{m+1}\cdot \left(\frac{S_{m+1}}{S_{m}} \right)^{j-m-1} \, . 
\]

This  implies that 
\begin{eqnarray*}
m &\le& \Exp(X) \\
&\le& \Exp\bigl(\max\{X\,,\, m-d-1\} \bigr)\\ &=& m-d-1 + \sum_{j\ge m-d} S_j \\
&\le& m-d-1 + \sum_{j= m-d}^{\infty}S_{m+1}\cdot\left(\frac{S_{m+1}}{S_{m}} \right)^{j-m-1} \\
&=& m-d-1 + \frac{S_{m+1}}{1 - \frac{S_{m+1}}{S_m}} \cdot \left(\frac{S_{m+1}}{S_m}\right)^{-d-1} \, .
\end{eqnarray*}
Summarizing the above, we  have shown that  
\[
S_{m+1} \ge (d+1) \left(1 - \frac{S_{m+1}}{S_m}\right)\cdot 
\left(\frac{S_{m+1}}{S_{m}} \right)^{d+1} \, .
\]
Since    $d\le \frac{S_{m+1}}{p_m} < d+1$,  it follows that $\frac{d}{d+1} \le \frac{S_{m+1}}{S_m}\le \frac{d+1}{d+2}$ which,  combined with the fact that the function $f(x)=(1-x)\cdot x^{d+1}$ is   increasing in the interval $[\frac{d}{d+1}, \frac{d+1}{d+2}]$, yields 
$ \left(1 - \frac{S_{m+1}}{S_m}\right)\cdot 
\left(\frac{S_{m+1}}{S_{m}} \right)^{d+1}\ge f(\frac{d}{d+1})$. The desired result follows upon observing that  
\[
S_{m+1} \ge (d+1) \cdot \left(1-\frac{d}{d+1}\right) \cdot \left(\frac{d}{d+1} \right)^{d+1} = \left(1-\frac{1}{d+1} \right)^{d+1} \ge \frac{1}{4} \, ,
\]
where the last inequality follows from the fact that the sequence $\{(1-1/r)^{r}\}_{r\ge 2}$ is increasing. 
Since $\gamma\le1/2$, it follows that 
\[
S_{m+1}\ge\frac14\ge\frac{\gamma}{1+2\gamma}\, , 
\]
and the first assertion of
Theorem~\ref{thm:main1} follows.

For the second assertion, set
$\mu=\Exp(X)$ and $m=\lfloor\mu\rfloor$.
Since an IFR distribution has interval support, it follows that 
$p_m>0$. Thus~\eqref{eq:left_skewed} implies that 
$p_{m+1}\ge p_m/2>0$ and the first assertion, with $\gamma=1/2$, yields 
\[
\Prb(X>\mu)=\Prb(X\ge m+1)\ge\frac14 \, .
\]
This completes the proof of Theorem~\ref{thm:main1}.

\section{Proof of Theorem~\ref{thm:main_hyp}}
\label{sec:hypergeometric}

\noindent
Let $H\sim\Hyp(n,i,k)$. 
To simplify notation, set $\mu_H=\Exp(H)$ and $\lambda_H=\frac{(n-i)(n-k)}{n} = n-i-k+\mu_H$. 
Moreover, set $p_j=\Prb(H=j)$ for $j\in [\max\{0,k-(n-i)\}\,..\,\min\{k,i\}]$.  
Hence we may rewrite the hypothesized inequalities as 
\begin{equation}\label{eq:check1}
\mu_H\,\ge\, \log(4/3)  \qquad\text{ and }\qquad \lambda_H \,\ge\, \log(4/3) \, .
\end{equation}
The hypotheses imply that $H$ is nondegenerate. Its support is
a discrete interval, and
\[
\frac{p_{j+1}}{p_j}=\frac{(i-j)(k-j)}{(j+1)(n-i-k+j+1)}
\]
is non-increasing in $j$ wherever the adjacent probabilities are positive.
It follows that $H$ is log-concave and therefore IFR.

We begin with the proof of first  assertion in Theorem~\ref{thm:main_hyp}; hence we assume that $\mu_H$ is not an integer. Assume first that  $0<\mu_H<1$. Since
$n-i-k$ is an integer such that 
$n-i-k=\lambda_H-\mu_H>-1$, it holds $n-i-k\ge0$. Hence 
\[
p_0=\prod_{j=0}^{k-1}\left(1-\frac{i}{n-j}\right)
\le e^{-ik/n}=e^{-\mu_H}\le\frac34,
\]
and $\Prb(H>\mu_H)=1-p_0\ge1/4$, and the desired result follows. 
Similarly, if $0<\lambda_H<1$ then, upon considering the random variable $H' = n-i-k+H\sim\Hyp(n,n-i,n-k)$, we have 
$\mu_{H'},\lambda_{H'}\ge \log(4/3)$ as well as 
$\mu_{H'}=\lambda_H$, and the already proven case implies that $\Prb(H> \mu_H) = \Prb(H' > \mu_{H'})\ge 1/4$.

Hence we may assume that $\mu_H,\lambda_H>1$.
Set $m=\lfloor\mu_H\rfloor$, and note that  $p_m,p_{m+1}>0$. 
Moreover, it holds $m\ge1$,\, $ik-nm\ge1$ and, since $\lambda_H>1$ and $0<\mu_H-m<1$,  
$n-i-k+m=\lambda_H-(\mu_H-m)>0$. We now compute 
\[
2(i-m)(k-m)-(m+1)(n-i-k+m+1)=(m-1)(n-i-k+m-1)+2(ik-nm-1) \ge 0 \, ,
\]
and thus $p_{m+1}/p_m\ge1/2$. Theorem~\ref{thm:main1},
with $\gamma=1/2$, now implies that 
\[
\Prb(H>\mu_H)=\Prb(H\ge m+1)\ge\frac14 \, ,
\]
as desired. 

We now prove the second assertion; so we assume that $\mu_H$ is an integer. Hence 
$\lambda_H=n-i-k+\mu_H$ is also an integer and  $\mu_H,\lambda_H\ge1$.
Since $H$ is non-degenerate with a discrete interval support, it holds that 
$p_{\mu_H},p_{\mu_H+1}>0$. Moreover, the identity $ik=n\mu_H$, implies 
\[
(i-\mu_H)(k-\mu_H)=\mu_H\cdot(n-i-k+\mu_H)=\mu_H\cdot\lambda_H \, .
\]
Hence, it holds 
\[
\frac{\Prb(H = \mu_H+1)}{\Prb(H=\mu_H) }
=\frac{\mu_H\cdot \lambda_H}{(\mu_H+1)(\lambda_H+1)}
\ge\frac{\mu_H}{2(\mu_H+1)} \, ,
\]
and Theorem~\ref{thm:main1}, with
$\gamma=\frac{\mu_H}{2(\mu_H+1)}\in(0,1/2)$, yields 
\[
\Prb(H>\mu_H)\ge\frac{\gamma}{1+2\gamma}=\frac{\mu_H}{4\mu_H+2}
\ge\frac16 \, . 
\]
This proves the second assertion.

It remains to show that the constants are best possible. 
To prove that the constant $1/6$ in the second assertion is best possible, 
for each integer $r\ge1$, consider the random variable 
$H_r\sim\Hyp(2r+2,2r,r+1)$. 
Then $\Exp(H_r)=r$,\, $\lambda_{H_r}=1$ and 
\[ 
\Prb(H_r>r) =\frac{\binom{2r}{r+1}}{\binom{2r+2}{r+1}} =\frac{r}{4r+2} \, .
\]
This shows that the bound $\frac{\mu_H}{4\mu_H+2}$ is best possible,
including the constant $1/6$ when $r=1$.

To prove that the constant $1/4$ in the first assertion is best possible, let $r\ge 3$ be an integer and consider a $H_r\sim\Hyp(2r, r+1, 2)$ random variable. 
Then $\frac{2(r+1)}{2r} = 1 + \frac{1}{r}>1$ and $\frac{(r-1)(2r-2)}{2r} = \frac{(r-1)^2}{r}>1$. Hence~\eqref{eq:check1} holds true, and 
\[
\Prb(H_r > \Exp(H_r)) = \Prb(H_r =2) = \frac{\binom{r+1}{2}}{\binom{2r}{2}}  =\frac{r+1}{2(2r-1)}\longrightarrow\frac14 \, , 
\]
as $r\to\infty$.

Finally, to show that the constant $\log(4/3)$ is best possible,  
fix $\delta \in (0, \log(4/3))$. For each integer  $r\ge 3$, set $i_r= \lceil\delta r^2\rceil$ and consider a random variable $H_r \sim \Hyp(r^3,i_r, r)$. Note that $\Exp(H_r) = \frac{i_r}{r^2}\to \delta$, as $r\to\infty$, and $\frac{(r^3 -i_r)(r^3-r)}{r^3}>1$. We will show that $\Prb(H_r >\Exp(H_r))< 1/4$ for sufficiently large $r$. 
To this end, note that 
for large enough $r$ it holds $\Exp(H_r)< \log(4/3)$, and so 
\[
\Prb\big(H_r \,>\, \Exp(H_r)\big) = \Prb(H_r\ge 1) = 1 - \Prb(H_r=0)  = 1 - \prod_{j=0}^{r-1} \left(1 - \frac{i_r}{r^3 - j} \right) \, .
\]
Now, using the inequality $\log(1-x)\ge -\frac{x}{1-x}$, for $x\in (0,1)$,   we deduce that  for large enough $r$ it holds 
\[
\log\big(\Prb(H_r=0)\big)\,=\, \sum_{j=0}^{r-1} \log\left(1 - \frac{i_r}{r^3 - j} \right) \,\ge\, - \sum_{j=0}^{r-1} \frac{i_r}{r^3-i_r-j} \ge - \frac{r\cdot i_r}{r^3 - i_r-r+1} \,\to\, -\delta \, ,
\]
as $r\to \infty$. Hence, for large enough $r$, it holds  
\[
\Prb(H_r\ge 1) = 1 - \Prb(H_r=0) \le 1 - e^{- \frac{r\cdot i_r}{r^3 - i_r-r+1}} \,\to\, 1 - e^{-\delta} < \frac{1}{4} \, ,  
\]
and the proof of Theorem~\ref{thm:main_hyp} is complete.

\section{Some applications} 
\label{sec:applications}

In this section we apply Theorem~\ref{thm:main1} to several
discrete distributions. We begin with sums of independent
Bernoulli random variables whose success probabilities
need not be the same. 

\begin{proposition}
\label{prop:poisson_binomial_quarter}
Let $B_1,\ldots,B_n$ be independent Bernoulli random
variables with success probabilities
$0<p_1\le p_2\le\cdots\le p_n<1$. Consider the sum  
$X=\sum_{i=1}^{n} B_i$ and set $\mu = \Exp(X)=\sum_{i=1}^{n}p_i$ and $m=\lfloor\mu\rfloor$. 
If either\, $\log(4/3)\le\mu<1$ or 
\begin{equation}
\label{eq:bernoulli_odds_condition}
\sum_{i=1}^{n-m}\frac{p_i}{1-p_i}
\ge \frac{m+1}{2} 
\end{equation}
hold true, then 
\[
\Prb(X>\Exp(X))\ge\frac14 \, .
\]
\end{proposition}
\begin{proof}
Note that the assumption that the sequence $\{p_i\}_{i=1}^{n}$ is non-decreasing implies that the sequence $\{\frac{p_i}{1-p_i}\}_{i=1}^{n}$ is non-decreasing. 

If $\log(4/3)\le\mu<1$, then it follows from the independence  of the $B_i, i\in [n]$, that 
\[
\Prb(X>\mu) = 1-\Prb(X=0) = 1-\prod_{i=1}^{n}(1-p_i) \ge 1-e^{-\mu} \ge \frac14 \, .
\]

Now suppose that~\eqref{eq:bernoulli_odds_condition} holds.
Since Bernoulli distributions are log-concave and
log-concavity is preserved under convolution, $X$ is
log-concave and hence IFR.
Moreover, it holds $0\le m\le n-1$, and $\Prb(X=m)>0$.
For every $i\in [n]$, the independence of $B_1,\ldots,B_n$ implies that
\[
\Prb(X=m+1\,,\, B_i=1) = \frac{p_i}{1-p_i}\cdot \Prb(X=m\, ,\, B_i=0) \, .
\]
Hence, upon summing over all $i=1,\ldots,n$, we deduce that 
\begin{eqnarray*}
(m+1)\cdot \Prb(X=m+1)
&=&
\sum_{i=1}^{n}\frac{p_i}{1-p_i}\cdot\Prb(X=m\, ,\, B_i=0)\\
&=&
\Exp\left(
\mathbf{1}_{\{X=m\}}
\sum_{i:\,B_i=0}\frac{p_i}{1-p_i}
\right)
\,\ge\,
\Prb(X=m) \sum_{i=1}^{n-m}\frac{p_i}{1-p_i} \, ,
\end{eqnarray*}
where the last inequality holds because, on the event $\{X=m\}$,
exactly $n-m$ of the indicators vanish, and the  sequence $\{\frac{p_i}{1-p_i}\}_{i=1}^{n}$ is non-decreasing.
Consequently,
\[
\frac{\Prb(X=m+1)}{\Prb(X=m)}
\ge
\frac{1}{m+1}
\sum_{i=1}^{n-m}\frac{p_i}{1-p_i}
\ge\frac12.
\]
Theorem~\ref{thm:main1}, applied for $m=\lfloor\mu\rfloor$,
now yields the desired result.
\end{proof}

\begin{remark}
In case $p_1=\cdots=p_n=p$ in Proposition~\ref{prop:poisson_binomial_quarter}, then $X\sim\Bin(n,p)$ and Proposition~\ref{prop:poisson_binomial_quarter} recovers the bound $\Prb(X> \Exp(X))\ge 1/4$, for $p\in [\log(4/3)/n, 1)$,  due to Pinelis~\cite{Pinelis}. 
Indeed, 
if $\log(4/3)\le np<1$, the result follows immediately from
Proposition~\ref{prop:poisson_binomial_quarter}.
If $np\ge1$, then set $m=\lfloor np\rfloor$ and note that, 
since all success probabilities equal $p$, it holds 
$\sum_{i=1}^{n-m}\frac{p}{1-p} = \frac{(n-m)p}{1-p} \ge m \ge\frac{m+1}{2}$. 
Thus~\eqref{eq:bernoulli_odds_condition} holds, and
Proposition~\ref{prop:poisson_binomial_quarter}
again gives the desired bound, which can also be deduced directly from Theorem~\ref{thm:main1}, using $\gamma=1/2$. 
\end{remark}

As another illustration of our main result, we 
apply it to the difference of independent
binomial random variables.

\begin{corollary}
\label{cor:binomial_difference}
Let $X_1\sim\Bin(n,p)$ and $X_2\sim\Bin(k,q)$ be independent,
where $n,k\ge1$ and $p,q\in(0,1)$.
Consider the difference $X = X_1-X_2$, and set $m=\lfloor np-kq\rfloor$. 
If
\begin{equation}
\label{eq:binomial_difference_condition}
2(n-m)\cdot \min\left\{\frac{p}{1-p},\frac{1-q}{q}\right\} \ge k+m+1 \, 
\end{equation}
holds true, then
\[
\Prb\bigl(X>\Exp(X)\bigr)\ge\frac14 \, .
\]
The constant $1/4$ is best possible.
\end{corollary}

\begin{proof}
The random variable
$Y=X+k=X_1+(k-X_2)$ is a sum of $n+k$ independent Bernoulli variables; $n$ having success probability $p$, and $k$ having success
probability $1-q$.
Moreover, it holds 
$\lfloor\Exp(Y)\rfloor=m+k$. Now observe that each ``success odds" is either equal to $\frac{p}{1-p}$ or equal to $\frac{1-q}{q}$. 
This implies that the sum of the 
$(n+k)-(m+k)=n-m$ smallest ``success odds" is at least
\[
(n-m)\cdot \min\left\{\frac{p}{1-p},\frac{1-q}{q}\right\} \ge\frac{k+m+1}{2} \, .
\]
Proposition~\ref{prop:poisson_binomial_quarter} applies to $Y$, and the observation that 
$\Prb(X>\Exp(X))=\Prb(Y>\Exp(Y))$ gives the desired bound. 

To see that the constant $1/4$ is best possible, take $n=k=1$ and $p=q=1/2$.
Then~\eqref{eq:binomial_difference_condition} holds with
equality, $\Exp(X)=0$, and
$\Prb(X>0) = \Prb(X_1=1,X_2=0) = \frac14$. 
\end{proof}

Further applications to various discrete distributions can be obtained by verifying condition~\eqref{eq:criterion-assumptions}. We illustrate some of those applications in the following.

\begin{corollary}\label{cor:further}
Let $c_r=r\bigl((4/3)^{1/r}-1\bigr)$ for real $r\ge1$; note that $c_1=1/3$, that $c_r$ decreases
to $\log(4/3)$, and that $(1+\mu/r)^{-r}\le 3/4$ if and only if $\mu\ge c_r$. The following hold true. 
\begin{enumerate}
\item[(a)] If $X\sim\Pois(\lambda)$ with $\lambda\ge\log(4/3)$, then $\Prb(X>\lambda)\ge\frac14$.
\item[(b)] If $X\sim\mathrm{NB}(r,p)$ is a negative binomial random variable with $r\ge1$ and $p\in(0,1)$, then
  $\Prb\bigl(X>r(1-p)/p\bigr)\ge\frac14$ whenever $p^r\le\frac34$, i.e.\ whenever $\Exp(X)\ge c_r$.
\item[(c)] Let  $X$ count the failures before the $r$-th success when sampling without replacement
  from an urn having $a$ successes and $b$ failures, where $a,b\ge1$ and $1\le r\le a$ are integers.
  Then $\Prb\bigl(X>rb/(a+1)\bigr)\ge\frac14$ whenever \, $rb/(a+1)\ge c_r$.
\end{enumerate}
\end{corollary}

\begin{proof} 
(a) The ratio $\frac{p_{j+1}}{p_j}=\frac{\lambda}{j+1}$ is decreasing in $j$, so $X$ is log-concave. If $\lambda<1$,
then $1-p_0=1-e^{-\lambda}\ge1/4$. If $\lambda\ge1$, then set $m=\lfloor\Exp(X)\rfloor\ge 1$ and note that  $\frac{p_{m+1}}{p_m}=\frac{\lambda}{m+1}\ge \frac{m}{m+1}\ge \frac 12$. The result follows from Theorem~\ref{thm:main1}.

(b) Here $\frac{p_{j+1}}{p_j}=\frac{(1-p)(j+r)}{j+1}$, which is non-increasing in $j$ because $r\ge1$. Set 
$\mu=\Exp(X)=\frac{r(1-p)}{p}$, so that $1-p=\frac{\mu}{r+\mu}$ and $p^r=(1+\mu/r)^{-r}$. If $\mu<1$, then
$1-p_0=1-p^r\ge1/4$. If $\mu\ge1$, then set $m=\lfloor\mu \rfloor\ge 1$ then and note that  
$\frac{p_{m+1}}{p_m}=\frac{\mu(m+r)}{(r+\mu)(m+1)}\ge\frac{m}{m+1}\ge\frac{1}{2}$, where the first inequality follows from the fact that it is equivalent  $\mu\ge m$. The desired result follows from Theorem~\ref{thm:main1}. 

(c) The probability mass function is equal to $p_j=\binom{j+r-1}{j}\binom{a+b-r-j}{a-r}\big/\binom{a+b}{a}$, for
$j\in[0..b]$. Also,  it holds $\mu=\Exp(X)=\frac{rb}{a+1}$.  Now observe that 
$\frac{p_{j+1}}{p_j}=\frac{j+r}{j+1}\cdot\frac{b-j}{a+b-r-j}$ is a product of two positive non-increasing
factors (the second because $a\ge r$);  hence $X$ is log-concave. 
If $\mu<1$, then the fact that $\frac ba = \frac{\mu(a+1)}{ra}>\frac{\mu}{r}$ implies that 
\[
p_0=\prod_{i=0}^{r-1}\frac{a-i}{a+b-i} = \prod_{i=0}^{r-1} \frac{1}{1+ \frac{b}{a-i}} \le\bigl(1+\tfrac  ba\bigr)^{-r}<\bigl(1+\tfrac{\mu}{r}\bigr)^{-r}\le\tfrac34 \, ,
\]
If $\mu\ge1$, then set $m=\lfloor\mu\rfloor$ and note that $m \le\mu  = b\cdot\frac{r}{a+1}< b$; hence $rb-ma\ge m$. Moreover, it holds   $(m+r)(b-m)-m(a+b-r-m)=rb-ma\ge m>0$, which yields   $\frac{p_{m+1}}{p_m}- \frac{m}{m+1}= \frac{rb-ma}{(m+1)(a+b-r-m)} > 0$. 
\end{proof}

\section*{Acknowledgements}
JA and OK were funded by the Czech Science Foundation project 24-11820S (``Automatic Combinatorialist").

\par\smallskip
\noindent
\textbf{AI Disclosure.}   
The authors used ChatGPT 5.6 Sol, ChatGPT 6.0 Astra, and Claude Fable 5.1 to assist in discovering preliminary proofs of the hypergeometric “quarter-bound.” These proofs were not obtained in one shot, but emerged through an iterative process in which the authors contributed mathematical ideas, evaluated the proposed arguments, and guided the AI tools from an initially cumbersome, computer-assisted proof toward an approach closer to the one presented in the manuscript. The authors subsequently simplified and generalized the results, making occasional further use of these tools for suggestions, verification, and literature searches. The manuscript itself was written by the human authors, who take full responsibility for its mathematical content.


\begin{thebibliography}{99}


 
\bibitem{Ai_Pelekis} J.~Ai, C.~Pelekis, \textit{On lower bounds for hypergeometric tails}, Statistics \& Probability Letters \textbf{236} (2026) Article 110798.  


\bibitem{Alqasem_et_al} A.~Alqasem, H.~Aravinda,  A.~Marsiglietti, J.~Melbourne, \textit{On a conjecture of Feige for discrete log-concave distributions}, SIAM Journal on Discrete Mathematics \textbf{38}(1) (2024) 93--102.

\bibitem{Barlow_Marshall} R.E.~Barlow, A.W.~Marshall, \textit{Bounds for Distributions with Monotone Hazard Rate, I}, The
Annals of Mathematical Statistics \textbf{35}(3) (1964) 1234--1257.

\bibitem{Barlow_Proschan} R.E.~Barlow, F.~Proschan, \textit{Mathematical Theory of Reliability}, Wiley 1965. 

\bibitem{Bentkus} V.~Bentkus, \textit{On Hoeffding’s inequalities}, Annals of Probability \textbf{32}(2) (2004) 1650--1673.






\bibitem{Grunbaum} B.~Gr\"unbaum, \textit{Partitions of mass-distributions and of convex bodies by hyperplanes},
Pacific J. Math. \textbf{10} (1960) 1257--1261.


\bibitem{Hui_Park} S.~Hui,  C.J.~Park, \textit{The Representation of Hypergeometric Random Variables using Independent Bernoulli Random Variables}, Communications in Statistics -- Theory and Methods, 43(19) (2014) 4103--4108.




\bibitem{Lovasz_Vempala} L.~Lov\'asz,  S.~Vempala, \textit{The geometry of logconcave functions and sampling algorithms}, Random Structures \& Algorithms \textbf{30} (2007) 307--358.



\bibitem{Mattner_Schulz} L.~Mattner, J.~Schulz, \textit{On normal approximations to symmetric hypergeometric laws}, Transactions of the American Mathematical Society \textbf{370}(1) (2018) 727--748.



\bibitem{Mattner} L.~Mattner, \textit{Teachable normal approximations to binomial and related probabilities or confidence bounds}, preprint (2025), arXiv:2503.20852v2 [stat.OT].



\bibitem{Pinelis} I.~Pinelis, \textit{Best lower bound on the probability of a binomial exceeding its expectation}, Statistics \& Probability Letters \textbf{179} (2021) 109224.

\bibitem{Shah} G.P.~Shah, \textit{Bounds for Lattice Distributions Having Monotone Hazard Rate with Applications}, 
Operations Research Center, University of California, Berkeley, Report ORC 65-12, June 1965.

\bibitem{Shaked_Shantikumar} M.~Shaked, G.J.~Shanthikumar, \textit{Stochastic Orders}, Springer, New York, 2007.

\bibitem{Shevtsova} I.G.~Shevtsova, \textit{On the absolute constants in the Berry-Esseen inequality and its structural
and nonuniform improvements}, (In Russian, with an English abstract.) Informatika i Ee Primeneniya
\textbf{7} (2013) 124–125.





\bibitem{Vatutin_Mikhailov} V.A.~Vatutin, V.G.~Mikhailov, \textit{Limit theorems for the number of empty cells in an
equiprobable scheme for group allocation of particles}, Theory Probab. Appl. \textbf{27} (1983) 734--743. Russian
original in Teor. Veroyatnost. i Primenen. \textbf{27} (1982) 684–692.




\end{thebibliography}
\end{document}